\documentclass[12pt,titlepage]{article}

\usepackage{graphicx,color}
\usepackage{amsmath,amssymb,amsthm}
\newtheorem{thm}{Theorem}
\newtheorem{lemma}[thm]{Lemma}
\newtheorem{cor}[thm]{Corollary}

\newtheorem{prop}[thm]{Proposition}

\newtheorem{defi}{Definition}

\def\beq{\begin{equation}}\def\eeq{\end{equation}}
\def\beqn{\begin{eqnarray}}\def\eeqn{\end{eqnarray}}

\def\qed{\ifhmode\unskip\nobreak\fi\quad\ifmmode\Box\else$\Box$\fi}

\newcommand\KG{{\rm KG}}

\title{On multichromatic numbers of widely colorable graphs}

\author{{\bf Anna Gujgiczer}$^{a,c}$
$\qquad$ {\bf G\'abor Simonyi}$^{b,a,}$
\\ \\
$^a$Department of Computer Science and Information Theory,\\
Faculty of Electrical Engineering and Informatics\\
Budapest University of Technology and Economics\\ \\
$^b$Alfr\'ed R\'enyi Institute of Mathematics, Budapest\\ \\
$^c$MTA-BME Lend\"ulet Arithmetic Combinatorics Research Group\\ \\
{\tt gujgicza@cs.bme.hu} \ \ \ {\tt simonyi@renyi.hu}}
\date{}

\begin{document}
\maketitle
\begin{abstract}
A coloring is called $s$-wide if no walk of length $2s-1$ connects vertices of the same color. A graph is $s$-widely colorable with $t$ colors if and only if it admits a homomorphism into a universal graph $W(s,t)$. Tardif observed that the value of the $r^{\rm th}$ multichromatic number $\chi_r(W(s,t))$ of these graphs is at least $t+2(r-1)$ and equality holds for $r=s=2$. He asked whether there is equality also for $r=s=3$. We show that $\chi_s(W(s,t))=t+2(s-1)$ for all $s$ thereby answering Tardif's question. We observe that for large $r$ (with respect to $s$ and $t$ fixed) we cannot have equality and that for $s$ fixed and $t$ going to infinity the fractional chromatic number of $W(s,t)$ also tends to infinity. The latter is a simple consequence of another result of Tardif on the fractional chromatic number of generalized Mycielski graphs.

\end{abstract}

\section{Introduction}
\message{Introduction}

A vertex-coloring of a graph is called $s$-wide if the two endvertices of every walk of length $2s-1$ receive different colors in it. If every vertex gets a different color then the coloring is $s$-wide if and only if the graph does not contain any odd cycle shorter than $2s+1$. The interesting phenomenon is that some graphs have $s$-wide colorings that are also optimal colorings.

A $1$-wide coloring is just a proper coloring. $2$-wide colorings were first investigated by Gy\'arf\'as, Jensen, and Stiebitz \cite{GyJS} who, answering a question of Harvey and Murty, showed the existence of a $t$-chromatic graph for every $t\ge 2$ with the additional property that it admits a $t$-coloring in which the neighborhood of every color class is an independent set. The analogous statement including more distant neighborhoods is also proven in \cite{GyJS}.

$3$-wide colorings (that are called simply wide colorings in \cite{ST06}) turned out to be relevant concerning the local chromatic number of several graph families whose chromatic number can be determined by the topological method of Lov\'asz \cite{LLKn}, cf. \cite{ST06} for more details and also for the relevance of $s$-wide colorability in the context of the circular chromatic number.

A graph homomorphism from a graph $F$ to a graph $G$ is an edge-preserving map of the vertex set of $F$ to the vertex set of $G$. The existence of such a map is denoted by $F\to G$. It is easy to see that $G\to K_t$ is equivalent to the $t$-colorability of graph $G$, that is to $\chi(G)\le t$, where $\chi(G)$ is the chromatic number of $G$. We refer to the book \cite{HN} for a general treatment of the theory of graph homomorphisms.

Several other types of graph colorings can also be expressed by the existence of a graph homomorphism to some target graph and $s$-wide colorability is no exception.
It is proven independently in \cite{BS} and \cite{ST06} (and already in \cite{GyJS} for the $s=2$ case) that $s$-wide colorability with $t$ colors is equivalent to the existence of a homomorphism to the following graph we denote by $W(s,t)$ as in \cite{ST06}.
$$V(W(s,t))=\{(x_1\dots x_t):\ \forall i\ x_i\in\{0,1,\dots,s\}, \exists! i\ x_i=0,\ \exists j\ x_j=1\},$$
$$E(W(s,t))=\{\{(x_1\dots x_t),(y_1\dots y_t)\}: \forall i\ |x_i-y_i|=1\ {\rm or}\ x_i=y_i=s\}.$$

\begin{prop} \label{prop:Wuniv}{\rm (\cite{BS, GyJS, ST06})}
A graph $G$ admits an $s$-wide coloring using $t$ colors if and only if $G\to W(s,t)$.
\end{prop}

\medskip

A different incarnation of the graphs $W(s,t)$ appears in the papers \cite{HH, Tardif05, MWjctb}, where (following Wrochna's notation in \cite{MWjctb}) a graph operation $\Omega_k$ is given for every odd integer $k$ and when applied to the complete graph $K_t$ for $k=2s-1$ then the resulting graph is isomorphic to $W(s,t)$. We will give and make use of this alternative definition in Section~\ref{sect:main}.

\medskip

It is easy to see that $W(s,t)$ can be properly colored with $t$ colors: set the color of vertex $(x_1\dots x_t)$ to be the unique $i$ for which $x_i=0$.
It is proven in \cite{BS, GyJS, ST06} (cf. also the chromatic properties of the more general $\Omega_k$ construction in \cite{HH, Tardif05, MWjctb})
that this coloring is optimal, that is,
\begin{equation}\label{eq:tchrom}
\chi(W(s,t))=t.
\end{equation}
This represents the surprising fact, that there are $t$-chromatic graphs that can be optimally colored in such a way, that the complete $d$-neighborhood of any color class is an independent set for every $d<s$. (By $d$-neighborhood of a color class we mean the set of vertices at distance exactly $d$ from the closest element of the color class. In fact, if $G$ is $s$-widely colored then not only the $d$-neighborhoods of color classes form independent sets for $d<s$ but all those vertices that can be attained via walks of length $d$ from the given color class.) The proof of $t$-chromaticity of $W(s,t)$ goes via showing that some other graphs that are known to be $t$-chromatic admit a homomorphism into $W(s,t)$. These graphs include generalized Mycielski graphs, Schrijver graphs, and Borsuk graphs of appropriate parameters (for the definition of generalized Mycielski graphs see Section~\ref{sect:elab}; cf. \cite{Sch78, EH, Mat} for the definition of Schrijver graphs and Borsuk graphs and \cite{ST06} for further details), showing in particular that all these graphs admit $s$-wide colorings. A common property of all these graphs is that their chromatic number can be determined by the already mentioned topological method introduced by Lov\'asz in his celebrated paper \cite{LLKn} proving Kneser's conjecture.

For $n,k$ positive integers satisfying $n\ge 2k$ the Kneser graph ${\rm KG}(n,k)$
is defined on ${\binom{[n]}{k}}$,
the set of all $k$-element subsets of the $n$-element set $[n]=\{1,2,\dots,n\}$ as vertex set. Two vertices are adjacent if and only
if the $k$-element subsets they represent are disjoint. It is not hard to show that $\chi({\rm KG}(n,k))\le n-2k+2$ (for all $n,k$ satisfying $n\ge 2k$) and Kneser \cite{Kne} conjectured that this estimate is sharp. This was proved by Lov\'asz \cite{LLKn} thereby establishing the following result.

\smallskip
\par\noindent
{\bf Theorem} (Lov\'asz-Kneser theorem). $$\chi({\rm KG}(n,k))= n-2k+2.$$

\smallskip
\par\noindent
For more about the topological method we refer to the excellent book by Matou\v{s}ek \cite{Mat}.

\medskip
\par\noindent
The existence of a homomorphism to the Kneser graph $\KG(n,k)$ can also be interpreted as a coloring property: $G\to\KG(n,k)$ holds if and only if we can color the vertices of $G$ with $n$ colors in such a way that every vertex receives $k$ distinct colors and if two vertices $u$ and $v$ are adjacent then the set of colors received by $u$ is disjoint from the set of colors received by $v$. Such colorings were first considered by Geller and Stahl, see \cite{GS, Stahl}. Stahl \cite{Stahl} introduced the corresponding chromatic number $\chi_k(G)$ as the minimum number of colors needed for such a coloring, called a $k$-fold coloring and $\chi_k(G)$ the $k$-fold chromatic number in \cite{SchU} (or $k$-tuple chromatic number in \cite{HN}). The fractional chromatic number $\chi_f(G)$ can be defined as $$\chi_f(G)=\inf_k\left\{\frac{\chi_k(G)}{k}\right\}=\inf\left\{\frac{n}{k}: G\to\KG(n,k)\right\}.$$
Note the immediate consequence of this definition that if $G\to H$ then $\chi_f(G)\le\chi_f(H)$.
\medskip
\par\noindent
Not surprisingly, determining multichromatic numbers (that is, $k$-fold chromatic numbers for various $k$'s) can be even harder in general than determining the chromatic number which is the special case for $k=1$.
An example of this phenomenon is that while the chromatic number of Kneser graphs is already known by the Lov\'asz-Kneser theorem, it is only a still open conjecture due to Stahl what homomorphisms exist and what do not between Kneser graphs, see Section 6.2 of \cite{HN} for details, cf. also \cite{TardZhu}.

\medskip
\par\noindent
The starting point of our investigations was a question by Tardif \cite{Tardnew} who observed that (\ref{eq:tchrom}) combined with the Lov\'asz-Kneser theorem implies that
\begin{equation}\label{eq:Tobs}
\chi_r(W(s,t))\ge t+2(r-1)
\end{equation}
 and that equality holds for $r=s=2$. (This is also true in the case of $r=s=1$ when it simply means $\chi(K_t)=t$.) Tardif asked if there is equality also for $r=s=3$. In particular, he was interested in whether
$W(3,8)\not\to {\rm KG}(12,3)$ and/or  $W(3,7)\not\to {\rm KG}(11,3)$ is true. Our main result will imply that this is actually not the case and equality does hold for $r=s=3$. The motivation for Tardif's question came from recent developments concerning Hedetniemi's conjecture in which wide colorings also turned out to be relevant.

\medskip
\par\noindent
Hedetniemi's conjecture asked whether the so-called categorical (or tensor) product $G\times H$ satisfies $\chi(G\times H)=\min\{\chi(G),\chi(H)\}.$
The conjecture is equivalent to say that $G\times H\to K_c$ implies that $G\to K_c$ or $H\to K_c$ must hold. (Although the latter directly only means $\chi(G\times H)\ge\min\{\chi(G),\chi(H)\}$, the reverse inequality is essentially trivial by $G\times H\to G$ and $G\times H\to H$ following easily from the definition of the categorical product.) If this holds for $K_c$, then $K_c$ is called multiplicative. Hedetniemi's conjecture  formulated in 1966 thus stated that $K_c$ is multiplicative for every positive integer $c$. This is trivial for $c=1$, easy for $c=2$ and is a far from trivial result by El-Zahar and Sauer \cite{EZS} for $c=3$ published in 1985. For no other $c$ it was decided (whether $K_c$ is multiplicative or not) until 2019, when a breakthrough by Yaroslav Shitov took place who proved in \cite{Shitov} that the conjecture is not true by constructing counterexamples for large enough $c$'s. The smallest $c$ for which Shitov's construction disproved the conjecture was extremely large (about $3^{95}$ according to an estimate in \cite{MWnew}). This value was dramatically improved within a relatively short time. Using Shitov's ideas in a clever way first Zhu \cite{Zhunew} reduced $c$ to $125$, then developing the method further Tardif \cite{Tardnew} showed a counterexample for $c=13$. He remarked that his construction would also work for $c=12$ and $11$, respectively, provided that $W(3,8)\not\to {\rm KG}(12,3)$ and  $W(3,7)\not\to {\rm KG}(11,3)$.
Our main result is the following that shows as a special case that these homomorphisms do exist.

\begin{thm}\label{thm:main}
$$\chi_s(W(s,t))=t+2(s-1).$$
\end{thm}

\par\noindent
Later Wrochna \cite{MWnew} managed to improve on Tardif's result using the ideas in \cite{Tardnew} in a different way and proving that $K_c$ is not multiplicative for any $c\ge 5$ thus leaving $c=4$ the only open case. (For more details about Hedetniemi's conjecture see e.g. Tardif's survey \cite{Tardsur} and the more recent papers cited above.)

\medskip
\par\noindent
The paper is organized as follows. We present the proof of Theorem~\ref{thm:main} in Section~\ref{sect:main}.
In Section~\ref{sect:elab} we elaborate on the problem of what we can say about $\chi_r(W(s,t))$ for general $r$. It will be an immediate consequence of Theorem~\ref{thm:main} combined with Tardif's observation (\ref{eq:Tobs}) that $\chi_r(W(s,t))=t+2(r-1)$ whenever $r\le s$. We will also observe
that we cannot have equality in (\ref{eq:Tobs}) for large enough $r$.
We will also show that the fractional chromatic number of $W(s,t)$ goes to infinity when $t$ grows and $s$ remains fixed. The paper is concluded with some observations about the position of the graphs $W(s,t)$ in the homomorphism order of graphs.

\section{Proof of the main result} \label{sect:main}

\medskip
\par\noindent
First we give the alternative definition of the graphs $W(s,t)$ using the graph operation $\Omega_k$, where $k=2\ell+1$ is odd, that was already mentioned in the Introduction. We give the definition of only $\Omega_{2\ell+1}(K_t)$ that we will use and refer to \cite{MWjctb} for the construction $\Omega_{2\ell+1}(G)$ for general graphs $G$.

\medskip
\par\noindent
\begin{defi}\label{defi:Omega}
The graph $\Omega_{2\ell+1}(K_t)$ is defined as follows.
%In the following definition all $A_i\subseteq [t]=\{1,\dots,t\}$.
$$V(\Omega_{2\ell+1}(K_t))=$$$$\{(A_0,A_1,\dots,A_{\ell}): \forall i\ A_i\subseteq [t], |A_0|=1, A_1\neq\emptyset, \forall i\in\{0,\dots,\ell-2\}\ A_i\subseteq A_{i+2}, A_{\ell-1}\cap A_{\ell}=\emptyset\},$$
%where $A_i\subseteq [t]=\{1,\dots,t\}$ for all $i\in\{0,1,\dots,\ell\}$.
\smallskip
$$E(\Omega_{2\ell+1}(K_t))=\{\{(A_0,A_1,\dots,A_{\ell}),(B_0,B_1,\dots,B_{\ell}):$$$$\ \forall i\in\{0,1,\dots,\ell-1\}\ A_i\subseteq B_{i+1}, B_i\subseteq A_{i+1}\ {\rm and}\ A_{\ell}\cap B_{\ell}=\emptyset\}.$$
\end{defi}
\medskip
\par\noindent
Note that the above conditions also imply that $A_{i-1}\cap A_i=\emptyset$ for all $1\le i\le \ell$ whenever $(A_0,A_1,\dots,A_{\ell})\in V(\Omega_{2\ell+1}(K_t))$.
\medskip
\par\noindent
It is straightforward and well-known (see e.g. \cite{MWjctb, MWnew}) that we have $$W(s,t)\cong \Omega_{2s-1}(K_t).$$
\smallskip
\par\noindent
%The above fact is straightforward and well-known (see e.g. \cite{MWjctb}), nevertheless we give the isomorphism (leaving it to the reader to see why it is an isomorphism) for the sake of completeness.
Indeed, one can easily check that the following function $g: V(W(s,t))\to V(\Omega_{2s-1}(K_t))$ provides an isomorphism between $W(s,t)$ and $\Omega_{2s-1}(K_t)$.
$$g: (x_1\dots x_t)\mapsto (A_0,A_1,\dots,A_{s-1}),$$
where
$$\forall i\in \{0,1,\dots,s-1\}: A_i=\{j: x_j\le i\ {\rm and}\ x_j\equiv i\ {\rm mod\ 2}\}.$$

\medskip
\par\noindent
{\it Remark 1.} We gave both descriptions of the graphs $W(s,t)$, because we believe that both are useful. In particular, we will formulate the proof of Theorem~\ref{thm:main} using the description of $\Omega_{2s-1}(K_t)$ as we believe that it makes the presentation of the proof easier to follow. Nevertheless, when we were thinking about the proof we felt we could understand the structure of these graphs better by considering its vertices as the sequences given in its definition as $W(s,t)$.
(It is also remarked in \cite{MWnew} that it is the $W(s,t)$ type description from which one easily sees that the number of vertices is $t(s^{t-1}-(s-1)^{t-1})$.)
$\Diamond$

\medskip
\par\noindent
Next we recall Tardif's observation (\ref{eq:Tobs}) that we state as a lemma for further reference and also prove for the sake of completeness.

\medskip
\par\noindent
\begin{lemma}\label{lem:Tobs} {\rm (Tardif \cite{Tardnew})}
For all positive integers $r$ and $s$ $$\chi_r(W(s,t))\ge t+2(r-1).$$
\end{lemma}

\proof
We cannot have $W(s,t)\to\KG(t+h,r)$ for $h<2(r-1)$ as $\chi(\KG(t+h,r))=t+h-2r+2$ by the Lov\'asz-Kneser theorem and this value is less then $t=\chi(W(s,t))$ whenever $h<2(r-1)$.
\hfill$\Box$

\bigskip
\par\noindent
{\bf Proof of Theorem~\ref{thm:main}.}
We need to show $$\chi_s(W(s,t))=\chi_s(\Omega_{2s-1}(K_t))=t+2(s-1).$$
Lemma~\ref{lem:Tobs} already shows that the right hand side is a lower bound thus our task is to prove the reverse inequality which is equivalent to the existence of a graph homomorphism from $W(s,t)\cong\Omega_{2s-1}(K_t)$ to ${\rm KG}(t+2(s-1),s).$ Below we give such a homomorphism $$f: (A_0,A_1,\dots,A_{s-1})\mapsto \{z_0,\dots,z_{s-1}\},$$ where
$\{z_0,\dots,z_{s-1}\}\in {\binom{[t+2(s-1)]}{s}}=V({\rm KG}(t+2(s-1),s))$.
To emphasize the mapping for $U=(A_0,A_1,\dots,A_{s-1})$ we will also use the notation $z_i=f_i(U)$ when $f((A_0,A_1,\dots,A_{s-1}))=\{z_0,\dots,z_{s-1}\}$. (Note that we do not assume that the $z_i$'s are monotonically increasing with respect to their indices, we only need that all of them are distinct for a given $f(U)=\{z_0,\dots,z_{s-1}\})$.

\medskip
\par\noindent
First assume that $s\ge 3$ is odd. (The $s=1$ case is a trivial special case of (\ref{eq:tchrom}).)

\smallskip
\par\noindent
For every even $i\in\{2,\dots, s-1\}$ we consider the three sets $A_{i-2}, A_{i-1}, A_i$ and for each such triple we define two elements of $f(U)$, namely $f_{i-1}(U)=z_{i-1}$ and $f_i(U)=z_i$ as follows.
According to the relative sizes of these three sets we will decide which of the elements $t+i-1, t+i, (t+s-1)+i-1=t+s+i-2$, and $(t+s-1)+i=t+s+i-1$ will be put into the set $f(U)$. For every even $i$ we will either put two of these elements into $f(U)$ or if not then we will find enough elements from $[t]$ to compensate this hiatus. This will give us $s-1$ distinct elements of $f(U)$. Finally we will define $f_0(U)$ as the missing $s$th element of $f(U)$. The rules are as follows.

\medskip
\par\noindent
i) If $|A_{i-2}|>|A_{i-1}|$ then let $f_{i-1}(U)=t+i-1$ and $f_i(U)=t+i$.
If $|A_{i-1}|>|A_i|$, then let $f_{i-1}(U)=t+s+i-2$ and $f_i(U)=t+s+i-1$. (Note that by $A_{i-2}\subseteq A_i$ at most one of the above two inequalities can hold so our definition is meaningful.)
\medskip
\par\noindent
ii) If $|A_{i-2}|<|A_{i-1}|<|A_i|$, then we must have $|A_i\setminus A_{i-2}|\ge 2$. In that case choose $2$ distinct elements of $A_i\setminus A_{i-2}$ (these will be elements from $[t]$) to be $f_{i-1}(U)$ and $f_i(U)$.

\medskip
\par\noindent
iii) If $|A_{i-2}|<|A_{i-1}|=|A_i|$, then $|A_i\setminus A_{i-2}|\ge 1.$ Let $f_{i-1}(U)$ be an arbitrary element of $A_i\setminus A_{i-2}$ and let
$$f_i(U)=\left\{\begin{array}{lll}t+s+i-2&&\hbox{if $\min(A_{i-1}\cup A_i)\in A_{i-1}$}
\\t+s+i-1&&\hbox{if $\min(A_{i-1}\cup A_i)\in A_i$.}
\end{array}\right.$$
Note that since $A_{i-1}\cap A_i=\emptyset$, $f_i(U)$ will be well defined.

\medskip
\par\noindent
iv) If $|A_{i-2}|=|A_{i-1}|<|A_i|$, then let
$$f_{i-1}(U)=\left\{\begin{array}{lll}t+i-1&&\hbox{if $\min(A_{i-2}\cup A_{i-1})\in A_{i-2}$}
\\t+i&&\hbox{if $\min(A_{i-2}\cup A_{i-1})\in A_{i-1}$.}
\end{array}\right.$$
Since $A_{i-2}\cap A_{i-1}=\emptyset$, $f_{i-1}(U)$ is well defined.
Let $f_i(U)$ be an arbitrary element of $A_i\setminus A_{i-2}$. Such a choice is possible as $A_i\setminus A_{i-2}\neq\emptyset$ in this case.

\medskip
\par\noindent
v) If $|A_{i-2}|=|A_{i-1}|=|A_i|$ (which means $A_i=A_{i-2}$) then let
$$f_{i-1}(U)=\left\{\begin{array}{lll}t+i-1&&\hbox{if $\min(A_{i-2}\cup A_{i-1})\in A_{i-2}$}
\\t+i&&\hbox{if $\min(A_{i-2}\cup A_{i-1})\in A_{i-1}$.}
\end{array}\right.$$
Let
$$f_i(U)=\left\{\begin{array}{lll}t+s+i-2&&\hbox{if $\min(A_{i-1}\cup A_i)\in A_{i-1}$}
\\t+s+i-1&&\hbox{if $\min(A_{i-1}\cup A_i)\in A_i$.}
\end{array}\right.$$

\medskip
\par\noindent
vi) Finally, let $f_0(U)$ be equal to the unique $h\in A_0.$
\medskip
\par\noindent
Note that by the above we have defined $f_j(U)$ for every $0\le j\le s-1$ and if $j\neq j'$ then $f_j(U)\neq f_{j'}(U)$ thus we have $f(U)\in V({\rm KG}(t+2(s-1),s)$ as needed. We have to prove that $f$ is indeed a graph homomorphism from $W(s,t)\cong\Omega_{2s-1}(K_t)$ to ${\rm KG}(t+2(s-1),s).$ We do this first and consider the case of even $s$ (that will be similar) afterwards.

\medskip
\par\noindent
Consider $U=(A_0,A_1,\dots,A_{s-1})$ and $U'=(B_0,B_1,\dots,B_{s-1})$. We have to show that if $f(U)\cap f(U')\neq\emptyset$, then $\{U,U'\}\notin E(\Omega_{2s-1}(K_t)).$
\medskip
\par\noindent
Assume that $f(U)\cap f(U')\neq\emptyset$ and we have $h\in f(U)\cap f(U')$ for some $h\in [t]$. Then we have $h$ appearing in some $A_j$ and some $B_k$, where both $j$ and $k$ are even. In particular, $h\in A_{s-1}\cap B_{s-1}$, thus $A_{s-1}\cap B_{s-1}\neq\emptyset$, therefore $U$ and $U'$ cannot be adjacent.

\medskip
\par\noindent
Now assume that $f(U)\cap f(U')\neq\emptyset$ but the intersection is disjoint from $[t]$ thus we have $t+d\in f(U)\cap f(U')$ for some $1\le d\le 2s-2$.

\medskip
\par\noindent
If $d$ is odd and $d\le s-1$, then $d=i-1$ for some even $2\le i\le s-1$, thus $t+d\in f(U)$ means $t+d=t+i-1=f_{i-1}(U)$. If this happens then either $|A_{i-2}|>|A_{i-1}|$ or $|A_{i-2}|=|A_{i-1}|$ and $\min(A_{i-2}\cup A_{i-1})\in A_{i-2}$. Similarly, $t+d=t+i-1\in f(U')$ implies that either $|B_{i-2}|>|B_{i-1}|$ or $|B_{i-2}|=|B_{i-1}|$ and $\min(B_{i-2}\cup B_{i-1})\in B_{i-2}$. Assume for contradiction that $\{U,U'\}$ is an edge of our graph $\Omega_{2s-1}(K_t)$. Then we must have $A_{i-2}\subseteq B_{i-1}$ and $B_{i-2}\subseteq A_{i-1}$ implying $$|A_{i-2}|\le |B_{i-1}|\le |B_{i-2}|\le |A_{i-1}|\le |A_{i-2}|,$$ therefore we must have equality everywhere. By $A_{i-2}\subseteq B_{i-1}$ and $B_{i-2}\subseteq A_{i-1}$ (that follows from $\{U,U'\}\in E(\Omega_{2s-1}(K_t))$) this implies $A_{i-2}=B_{i-1}$ and $B_{i-2}=A_{i-1}$ and therefore $j:=\min(A_{i-2}\cup A_{i-1})=\min(B_{i-2}\cup B_{i-1}).$ Our assumption on $d$ then implies both  $j\in A_{i-2}$ and $j\in B_{i-2}=A_{i-1}$ which is impossible by $A_{i-2}\cap A_{i-1}=\emptyset.$

\medskip
\par\noindent
The situation is similar for the other possible values of $d$. If $d=i\le s-1$ is even, then $t+d=t+i\in f(U)\cap f(U')$ for some adjacent vertices $U, U'$ would again imply $$|A_{i-2}|=|B_{i-1}|=|B_{i-2}|=|A_{i-1}|$$
and thus $A_{i-2}=B_{i-1}, B_{i-2}=A_{i-1}$ as above. Our assumption on $d$ now would imply for $j=\min(A_{i-2}\cup A_{i-1})=\min(B_{i-2}\cup B_{i-1})$ that it must be both in $A_{i-1}$ and in $B_{i-1}=A_{i-2}$ leading to the same contradiction as in the previous paragraph.
\smallskip
\par\noindent
For $s-1<d$ and $t+d\in f(U)\cap f(U')$ for adjacent vertices $U,U'$ we get the same contradiction with the indices shifted by one. In particular, this assumption implies $|A_{i-1}|\ge |A_i|$ and $|B_{i-1}|\ge |B_i|$ that by the adjacency of $U$ and $U'$ (meaning, in particular, $A_{i-1}\subseteq B_i$ and $B_{i-1}\subseteq A_i$) would imply
$$|A_{i-1}|=|B_i|=|B_{i-1}|=|A_i|$$ and thus $A_{i-1}=B_i$ and $B_{i-1}=A_i$. Then we obtain that $k:=\min(A_{i-1}\cup A_i)=\min(B_i\cup B_{i-1})$ should belong (depending on the parity of $d$) to both $A_{i-1}$ and $B_{i-1}=A_i$ or to both $A_i$ and $B_i=A_{i-1}$ leading to the same contradiction that $A_{i-1}\cap A_i\neq\emptyset$.
This finishes the proof for odd $s$.

\bigskip
\par\noindent
Now assume that $s$ is even. We need only some minor modifications compared to the odd $s$ case.
Let us now for every odd $i\in\{3,\dots,s-1\}$ define $f_{i-1}(U)$ and $f_i(U)$ almost the same way as in points i)-- v) above. (The only difference will be that the values $t+i-1$ and $t+i$ are shifted by $1$ to become $t+i$ and $t+i+1$. In case of $s=2$ the modified rules (i')-(v') will not apply, only those will that we denote by (vi') and (vii') below.) This gives the last $s-2$ values of the set $f(U)=\{f_0(U),f_1(U),\dots,f_{s-1}(U)\}$, what is left is to define $f_0(U)$ and $f_1(U)$ by a modified version of the sixth point above that has now two parts. The modified rules are as follows.

\medskip
\par\noindent
i') If $|A_{i-2}|>|A_{i-1}|$ then let $f_{i-1}(U)=t+i$ and $f_i(U)=t+i+1$.
If $|A_{i-1}|>|A_i|$, then let $f_{i-1}(U)=t+s+i-2$ and $f_i(U)=t+s+i-1$.
\medskip
\par\noindent
ii') and iii') are identical to ii) and iii), respectively.

\medskip
\par\noindent
iv') If $|A_{i-2}|=|A_{i-1}|<|A_i|$, then let
$$f_{i-1}(U)=\left\{\begin{array}{lll}t+i&&\hbox{if $\min(A_{i-2}\cup A_{i-1})\in A_{i-2}$}
\\t+i+1&&\hbox{if $\min(A_{i-2}\cup A_{i-1})\in A_{i-1}$.}
\end{array}\right.$$
%Since $A_{i-2}\cap A_{i-1}=\emptyset$, $f_{i-1}(U)$ is well defined.
Let $f_i(U)$ be an arbitrary element of $A_i\setminus A_{i-2}$.
%Such a choice is possible as $A_i\setminus A_{i-2}\neq\emptyset$ in this %case.

\medskip
\par\noindent
v)' If $|A_{i-2}|=|A_{i-1}|=|A_i|$ then let
$$f_{i-1}(U)=\left\{\begin{array}{lll}t+i&&\hbox{if $\min(A_{i-2}\cup A_{i-1})\in A_{i-2}$}
\\t+i+1&&\hbox{if $\min(A_{i-2}\cup A_{i-1})\in A_{i-1}$.}
\end{array}\right.$$
Let
$$f_i(U)=\left\{\begin{array}{lll}t+s+i-2&&\hbox{if $\min(A_{i-1}\cup A_i)\in A_{i-1}$}
\\t+s+i-1&&\hbox{if $\min(A_{i-1}\cup A_i)\in A_i$.}
\end{array}\right.$$

\medskip
\par\noindent
vi') If $|A_1|=|A_0|$, then let
$$f_0(U)=\left\{\begin{array}{lll}t+1&&\hbox{if $\min(A_0\cup A_1)\in A_0$}
\\t+2&&\hbox{if $\min(A_0\cup A_1)\in A_1$.}
\end{array}\right.$$
Note that in this case both $A_0$ and $A_1$ contains only one element and the value of $f_0(U)$ is $t+1$ or $t+2$ depending on which of the two is smaller. At the same time let $$f_1(U)=h\ {\rm where}\ A_1=\{h\},$$
that is, $h\in [t]$ is the unique element of $A_1$.

\medskip
\par\noindent
vii') If $|A_1|>|A_0|$, then since $|A_0|=1$ we have $|A_1|\ge 2$. Now choose two arbitrary distinct elements of $A_1$ for $f_0(U)$ and $f_1(U)$.
\medskip
\par\noindent
Note that we have $|A_1|\ge 1=|A_0|$ by the definition of $\Omega_{2s-1}(K_t)$, so we do not have to consider the possibility that $|A_0|>|A_1|$, it never occurs.

\medskip
\par\noindent
With this definition of $f(U)$ the proof that $f$ is a graph homomorphism is essentially identical to that we presented in the odd $s$ case. The main difference is that now those $j\in [t]$ that appear as elements of the sets $f(U)$ are all elements of some $A_i$ where $i$ is odd, while the corresponding $i$'s were all even in the case of odd $s$. The rest of the arguments work the same way as in the case of odd $s$.
\medskip
\par\noindent
This completes the proof.
\hfill$\Box$

\medskip
\par\noindent
We remark that by the composition of homomorphisms Theorem~\ref{thm:main} determines the $s$-fold chromatic number of every $s$-widely colorable $t$-chromatic graph.

\section{On other multichromatic numbers of $W(s,t)$}\label{sect:elab}

\medskip
\par\noindent
An immediate consequence of Theorem~\ref{thm:main} is that we can give the multichromatic numbers $\chi_r(W(s,t))$ for all $r\le s$.

\medskip
\par\noindent
\begin{cor}\label{cor:r<s}
If $r\le s$, then $$\chi_r(W(s,t))=t+2(r-1).$$
\end{cor}

\medskip
\par\noindent
The proof follows from the following simple lemma (which is essentially Lemma 2.3.(iv) of \cite{MWjctb}) combined with Tardif's observation given in Lemma~\ref{lem:Tobs}.

\medskip
\par\noindent
\begin{lemma}\label{lem:kisr}{\rm (\cite{MWjctb})}
For all $1\le r\le s$ we have $$W(s,t)\to W(r,t)$$
\end{lemma}

\proof
Define the following function for all $0\le a\le s$.
$$\varphi(a)=\left\{\begin{array}{lll}a&&\hbox{if $0\le a\le r$}
\\r&&\hbox{if $r<a\le s$.}
\end{array}\right.$$
It is straightforward to check that the mapping $g: (x_1\dots x_t)\mapsto (\varphi(x_1)\dots \varphi(x_t))$ is a homomorphism from $W(s,t)$ to $W(r,t)$ for all $1\le r\le s$.
\hfill$\Box$

\medskip
\par\noindent
{\bf Proof of Corollary~\ref{cor:r<s}}
In view of Lemma~\ref{lem:Tobs} it is enough to prove that $\chi_r(W(s,t))$ is at most the claimed value if $r\le s$. Applying Lemma~\ref{lem:kisr} and Theorem~\ref{thm:main} to $r\le s$ we have $$W(s,t)\to W(r,t)\to\KG(t+2(r-1),r)$$ implying $$\chi_r(W(s,t))\le t+2(r-1)$$ as needed.
\hfill$\Box$

\medskip
\par\noindent
For $r>s$ we do not know the value of $\chi_r(W(s,t))$. We know from Lemma~\ref{lem:Tobs} though that $\chi_r(W(s,t))\ge t+2(r-1)$ so the question naturally arises whether we could have equality here for every $r$. Below we show that this is not the case.
\medskip
\par\noindent
\begin{prop}\label{prop:nagyobb}
For all pairs of positive integers $t\ge 3$ and $s\ge 1$ there exists some threshold $r_0=r_0(s,t)>s$ for which
\begin{equation}\label{eq:nagyobb}
\chi_r(W(s,t))>t+2(r-1)
\end{equation}
whenever $r\ge r_0$.
\end{prop}

\proof
Assume for the sake of contradiction that for some fixed $s$ and $t$ we have $\chi_r(W(s,t))=t+2(r-1)$ for arbitrarily large $r$. That would imply that $\chi_f(W(s,t))\le\lim_{r\to\infty}\frac{t+2(r-1)}{r}=2.$
However, this cannot be true since $W(s,t)$ is not bipartite for $t\ge 3$ and thus it contains an odd cycle $C_{2b+1}$ for some positive integer $b$. Thus we must have $\chi_f(W(s,t))\ge\chi_f(C_{2b+1})=\frac{2b+1}{b}$, a number larger than $2$ with the constant value $\frac{1}{b}$.
\hfill$\Box$

\medskip
\par\noindent
The problem of determining the smallest possible $r$ for which (\ref{eq:nagyobb}) holds is left as an open problem. It is frustrating that we were not able to decide even whether this value is just $s+1$ as the proof of Theorem~\ref{thm:main} might suggest or larger.

\medskip
\par\noindent
{\it Remark 2.}
The previous proof does not specify $b$ as its value is not essential there. Nevertheless one can easily see that $W(s,3)\cong C_{6s-3}$. It is also easy to see that the odd girth of $W(s,t)$ must be at least $2s+1$ and we have equality here for $t\ge 2s+1$ since a cycle $C_{2s+1}$ is formed in $W(s,2s+1)$ by the vertices given by the sequence $(0,1,2,\dots,s,s,s-1,\dots,2,1)$ and its cyclic permutations. (For larger $t$ these sequences can be extended by an arbitrary number of coordinates equal to $s$.) In fact, the unpublished paper by Baum and Stiebitz \cite{BS} gives the general formula $2s-1+2\left\lceil\frac{2s-1}{t-2}\right\rceil$ for the odd girth of $W(s,t)$.
\hfill$\Diamond$

\medskip
\par\noindent
The previous proof raises the question what we can say about the fractional chromatic number of the graphs $W(s,t)$.
As a consequence of Theorem~\ref{thm:main} we know $\chi_f(W(s,t))\le\frac{t+2(s-1)}{s}$ and the previous simple proof implies that it is at least $2+\frac{1}{3s-2}$ for $t\ge 3$. Unfortunately we were not able to prove matching lower and upper bounds. But we can at least show that for any fixed $s$ the fractional chromatic number of $W(s,t)$ gets arbitrarily large as $t$ tends to infinity.

\medskip
\par\noindent
\begin{thm}\label{thm:megno}
For any fixed positive integer $s$ we have $$\lim_{t\to\infty}\chi_f(W(s,t))=\infty.$$
\end{thm}

\medskip
\par\noindent
The proof will be a simple consequence of the (already known) fact that certain generalized Mycielski graphs admit $s$-wide colorings.
To give more details we introduce generalized Mycielski graphs below.
\medskip
\par\noindent
\begin{defi}
The $h$-level generalized Mycielskian $M_h(G)$ of a graph $G$ is defined as follows.
$$V(M_h(G))=\{(v,j): v\in V(G), 0\le j\le h-1\}\cup\{z\}.$$
$$E(M_h(G))=\{\{(u,i),(v,j)\}: uv\in E(G)\ {\rm and}\ (|i-j|=1\ {\rm or}\ i=j=0\}\cup \{\{z,(v,(h-1))\}.$$
The $d$ times iterated $h$-level generalized Mycielskian $M_h(M_h(\dots M_h(G)\dots ))$ of a graph $G$ will be denoted by $M_h^{(d)}(G).$
\end{defi}

\medskip
\par\noindent
The term Mycielskian of a graph $G$ usually refers to $M(G)=M_2(G)$ and Mycielski graphs are the iterated Mycielkians of $K_2$ introduced by Mycielski \cite{Myc} as triangle-free graphs whose chromatic number grows by one at every iteration. The property $\chi(M(G))=\chi(G)+1$ is well-known to hold for any $G$ but the analogous equality is not always true for  $h$-level Mycielskians if $h>2$, cf. Tardif \cite{TardMyc}. Nevertheless Stiebitz \cite{Stieb} showed that $\chi(M_h(G))=\chi(G)+1$ is also true if $G$ is a complete graph or an odd cycle. (More generally one can say that this is the case whenever $G$ is a graph for which the topological lower bound on the chromatic number by Lov\'asz \cite{LLKn} is sharp, cf. \cite{GyJS, Mat} or \cite{ST06} for more details.)
 So by Stiebitz's result we have $$\chi(M_h^{(d)}(K_2))=d+2$$ for all positive integers $d$ and $h$.
\medskip
\par\noindent
The $t$-chromaticity of $W(s,t)$ is proven in \cite{BS, GyJS, ST06} by showing the existence of $t$-chromatic graphs that admit a homomorphism into $W(s,t)$. In case of \cite{BS, GyJS} these are generalized Mycielski graphs $M_h^{(t-2)}(K_2)$ for appropriately large $h$. (Since \cite{BS} is unpublished and \cite{GyJS} gives this explicitly only for $s=2$, we give some more details for the sake of completeness. Nevertheless, this is a straightforward generalization of the construction given in \cite{GyJS} as already noted in \cite{ST06} where the case $s=3$ is made explicit. So the following is a straightforward extension of Lemma~4.3 from \cite{ST06} also attributed to \cite{GyJS} there.)
\medskip
\par\noindent
\begin{lemma}\label{lem:gMycszeles} {\rm (\cite{GyJS})}
If $G$ has an $s$-wide coloring with $t$ colors, then $M_{3s-2}(G)$ has an $s$-wide coloring with $t+1$ colors.
\end{lemma}

\proof
Fix an $s$-wide coloring $c_0: V(G)\to [t]$ of $G$. Let $c: V(M_{3s-2}(G))\to [t]\cup\{\gamma\}$ be the following coloring using the additional color $\gamma$. Set $c(z)=\gamma$ and
$$c((v,j))=\left\{\begin{array}{lll}\gamma&&\hbox{if $j\in \{s,s+2,\dots,3s-4 \}$}\\
c_0(v)&&\hbox{otherwise.}\end{array}\right.$$
If we have a walk of odd length between vertices $(u,i)$ and $(v,j)$ with $c(u,i)=c(v,j)\in [t]$ that walk must either traverse the vertex $z$ or use an edge of the form $\{a,0),(b,0)\}$. In the latter case the walk
projects down to a walk of the same length between $u$ and $v$ in $G$ with $c_0(u)=c_0(v)$ so its length must be at least $2s+1$ by $c_0$ being $s$-wide. In case the walk traverses $z$ we can assume that we have $i\not\equiv j\ {\rm mod}\ 2$ and thus without loss of generality $j\equiv s\ {\rm mod}\ 2$ implying that $j\le s-2$. But then the distance between $(v,j)$ and $z$ is already at least $2s$, so the length of our walk is at least $2s+1$.
\smallskip
\par\noindent
Since deleting the set of vertices $\{(v,0)\}_{v\in V(G)}$ from $M_{3s-2}(G)$ the remaining induced subgraph is bipartite and $\gamma$ appears only on one side of this bipartite graph, any odd length walk between two vertices colored $\gamma$ must use an edge of the form $\{(u,0),(v,0)\}$. But the distance of any $\gamma$-colored vertex from such vertices is at least $s$, so such a walk also cannot be shorter than $2s+1$. Thus $c$ is indeed an $s$-wide coloring.
\hfill$\Box$

\medskip
\par\noindent
For $M(G)=M_2(G)$ Larsen, Propp and Ullman \cite{LPU} made the very nice observation, that $\chi_f(M(G))$ can be given by a simple function of $\chi_f(G)$, namely $$\chi_f(M(G))=\chi_f(G)+\frac{1}{\chi_f(G)}.$$ This was later generalized by Tardif for generalized Mycielskians.

\begin{thm}\label{thm:gMfrac} {\rm (Tardif \cite{TardMyc})}
$$\chi_f(M_h(G))=\chi_f(G)+{\frac{1}{\sum_{i=0}^{h-1}(\chi_f(G)-1)^i}}.$$
\end{thm}

\medskip
\par\noindent
Note that Tardif's theorem implies that $\chi_f(M_h^{(d)}(G))$ tends to infinity as $d$ goes to infinity for any fixed finite $h$.

\medskip
\par\noindent
{\bf Proof of Theorem~\ref{thm:megno}}
The proof is already immediate by the foregoing.
Lemma~\ref{lem:gMycszeles} and Tardif's Theorem~\ref{thm:gMfrac} together imply that $$\chi_f(W(s,t+1))\ge\chi_f(M_{3s-2}(W(s,t)))
=\chi_f(W(s,t))+{\frac{1}{\sum_{i=0}^{3s-3}(\chi_f(W(s,t))-1)^i}},$$
and this implies the statement.
\hfill$\Box$

\medskip
\par\noindent
In view of Lemma~\ref{lem:gMycszeles} it may be interesting to note that while a generalized Mycielskian of $W(s,t)$ admits a homomorphism into $W(s,t+1)$, the latter also admits a (very natural) homomorphism into another generalized Mycielksian of $W(s,t)$.

\medskip
\par\noindent
\begin{prop}\label{prop:masikhom}
$$W(s,t+1)\to M_s(W(s,t)).$$
\end{prop}

\proof
We explicitly give the homomorphism.
Let
$$g((x_1\dots x_{t+1}))=\left\{\begin{array}{lll} ((x_1\dots x_t),s-x_{t+1})&&\hbox{if $x_{t+1}>0\ {\rm and}\ (x_1\dots x_t)\in V(W(s,t))$ }\\
((01\dots 1),s-1)&&\hbox{if $\{i:x_i=1\}=\{t+1\}$ }\\
z&&\hbox{if $x_{t+1}=0$.}\end{array}\right.$$
(In fact, in the second case $((01\dots 1),s-1)$ can be substituted by an arbitrarily chosen $((y_1\dots y_t),s-1)$ for which $(y_1\dots y_t)\in V(W(s,t))$.)
\smallskip
\par\noindent
It is straightforward to check that the given function is indeed a graph homomorphism.
\hfill$\Box$

\medskip
\par\noindent
Thus we obtained that in the homomorphism order of graphs (cf. \cite{HN}) in which $F \preceq G$ if and only if $F\to G$ we have $W(s,t+1)$ sandwiched between two different generalized Mycialskians of $W(s,t)$, in particular,  $$M_{3s-2}(W(s,t))\preceq W(s,t+1)\preceq M_s(W(s,t)).$$
This excludes the possibility that our upper bound $\frac{t+2(s-1)}{s}$ on $\chi_f(W(s,t))$ provided by Theorem~\ref{thm:main} would be tight at least for all sufficiently large $t$, because then the difference $\chi_f(W(s,t+1))-\chi_f(W(s,t))$ would be equal to $\frac{1}{s}$ for large $t$ contradicting Tardif's Theorem~\ref{thm:gMfrac}.

\medskip
\par\noindent
With a little more considerations we can also show that $W(s,t+1)$ is actually {\em strictly} sandwiched between the above two generalized Mycielskians of $W(s,t)$ if $s>1$ and $t>2$.

\begin{prop} \label{prop:strict}
If $s\ge 2, t\ge 3$ then
\begin{equation}\label{eq:strict}
M_{3s-2}(W(s,t))\prec W(s,t+1)\prec M_s(W(s,t)).
\end{equation}
For $s=1$ all three graphs are isomorphic to $K_{t+1}$.
For $s>1, t=2$ we have
$$M_{3s-2}(W(s,2))\cong C_{6s-3}\cong W(s,3)\prec M_s(W(s,2))\cong C_{2s+1}.$$
\end{prop}

\proof
It is well-known and easy to prove that if $G$ is a vertex-color-critical graph (that is, one from which deleting any vertex its chromatic number decreases) and $\chi(M_h(G))=\chi(G)+1$, then $M_h(G)$ is also vertex-color-critical (see this e.g. as Problem 9.18 in the book \cite{LLbook} for $h=2$). It is shown independently both in \cite{BS} and \cite{ ST06} that $W(s,t)$ is edge-color-critical for every $s\ge 1,t\ge 2$.
Thus all three graphs appearing in (\ref{eq:strict}) are vertex-color-critical. Since they all have the same chromatic number this implies that any homomorphism that exists between any two of them should be onto. This also means that if any two of them would be homomorphically equivalent, then those two should have the same number of vertices, in particular, any homomorphism between them is a one-to-one mapping between their vertex sets. This is clearly not the case for the homomorphism given in the proof of Proposition~\ref{prop:masikhom} since several distinct vertices (their exact number is $s^t-(s-1)^t$) are mapped to the vertex $z$ unless $s=1$.
\smallskip
\par\noindent
If a homomorphism between $M_{3s-2}(W(s,t))$ and $W(s,t+1)$ was one-to-one then by the edge-color-criticality of $W(s,t+1)$ it cannot happen that we map two non-adjacent vertices of $M_{3s-2}(W(s,t))$ to two adjacent ones of $W(s,t+1)$, since then deleting the latter adjacency we would still have a homomorphism but into a graph of smaller chromatic number. Thus such a homomorphism would then be an isomorphism, that is the two graphs would be isomorphic which is clearly not the case if $s>1$ and $t>2$. (A quick way to see this is the following.
%The degree of vertex $(x_1\dots x__{t+1})$ is $k2^{t-k}$ where $k$ is the number of $1$'s among the $x_i$'s.
The maximum degree of $W(s,t+1)$ is $2^{t-1}$ attained by vertices $(x_1\dots x_{t+1})$ for which $|\{i: x_i=1\}|$ is equal to $1$ or $2$.
The maximum degree of $M_{3s-2}(W(s,t))$ is $|V(W(s,t))|=t(s^{t-1}-(s-1)^{t-1})$ that cannot be a power of $2$ for $s>1$ unless $t=2$.)
The remaining cases in the statement are straightforward to check.
\hfill$\Box$

\section{Acknowledgements}

\medskip
\par\noindent
We thank Claude Tardif for sharing with us already an early version of his paper \cite{Tardnew} which contained his interesting question that became the starting point of our work presented here.
\smallskip
\par\noindent
This research was partially supported by the National Research, Development and Innovation Office (NKFIH) grants K--120706 and BME NC TKP2020 and also by the BME-
Artificial Intelligence FIKP grant of EMMI (BME FIKP-MI/SC). The work of GS was also partially supported by the grants K--132696 and SSN-135643 of NKFIH Hungary.

\end{document}